\documentclass[conference]{IEEEtran}

\usepackage{cite}
\usepackage{amsmath,amssymb,amsfonts}
\usepackage{algorithmic}
\usepackage{graphicx}
\usepackage{textcomp}
\usepackage{xcolor}

\ifCLASSOPTIONcompsoc
    \usepackage[caption=false, font=normalsize, labelfont=sf, textfont=sf]{subfig}
\else
\usepackage[caption=false, font=footnotesize]{subfig}
\fi

\usepackage{mathtools} 
\usepackage{circuitikz} 

\usepackage{nomencl}
\makenomenclature

\usepackage{lipsum}

\newcommand{\refig}[1]{Fig.~\ref{#1}}
\newcommand{\retab}[1]{Table~\ref{#1}}

\newcommand{\e}[1]{\times10^{#1}}

\newcommand{\etal}{\textit{et al.}}

\newcommand{\defi}{\coloneqq}

\newcommand{\parentheses}[1]{\left( #1 \right)}

\newcommand{\iv}{I--V}

\newcommand{\Vpv}{V_\mathrm{pv}}
\newcommand{\Ipv}{I_\mathrm{pv}}

\newcommand{\Voc}{V_\mathrm{oc}}
\newcommand{\Isc}{I_\mathrm{sc}}
\newcommand{\Vmp}{V_\mathrm{mp}}
\newcommand{\Imp}{I_\mathrm{mp}}

\newcommand{\Iph}{I_\mathrm{ph}}
\newcommand{\Io}{I_\mathrm{o}}
\newcommand{\Rs}{R_\mathrm{s}}
\newcommand{\Rshh}{R_\mathrm{sh}}
\newcommand{\Gsh}{G_\mathrm{sh}}
\newcommand{\Ns}{N_\mathrm{s}}

\newcommand{\Id}{I_\mathrm{d}}
\newcommand{\Ish}{I_\mathrm{sh}}
\newcommand{\Vj}{V_\mathrm{j}}
\newcommand{\Tpv}{T_\mathrm{pv}}

\newcommand{\fsh}{f_\mathrm{sh}}
\newcommand{\fmp}{f_\mathrm{mp}}

\newcommand{\Rssh}{R_\mathrm{s}^\mathrm{sh}}
\newcommand{\Rsmp}{R_\mathrm{s}^\mathrm{mp}}

\newcommand{\Amax}{A_\mathrm{max}}
\newcommand{\Rsmin}{R_\mathrm{s,\,min}}

\newcommand{\Rsshlow}{R_\mathrm{s,\,low}^\mathrm{sh}}
\newcommand{\Rsmplow}{R_\mathrm{s,\,low}^\mathrm{mp}}

\newcommand{\Rsshhigh}{R_\mathrm{s,\,high}^\mathrm{sh}}
\newcommand{\Rsmphigh}{R_\mathrm{s,\,high}^\mathrm{mp}}

\newcommand{\ej}{ \exp{ \parentheses{ \cfrac{\Vj}{A} } } }

\newcommand{\eoc}{ \exp{ \parentheses{ \cfrac{\Voc}{A} } } }
\newcommand{\esc}{ \exp{ \parentheses{ \cfrac{\Rs\,\Imp}{A} } } }
\newcommand{\emp}{ \exp{ \parentheses{ \cfrac{\Vmp+\Rs\,\Imp}{A} } } }

\newcommand{\kb}{k_\mathrm{B}}

\usepackage{etoolbox}
\renewcommand\nomgroup[1]{%
  \item[\bfseries
  \ifstrequal{#1}{A}{Acronyms}{%
  \ifstrequal{#1}{P}{Physics constants}{%
  \ifstrequal{#1}{N}{Number sets}{%
  \ifstrequal{#1}{E}{Variables}{%
  \ifstrequal{#1}{S}{Subscripts}{%
  }}}}}%
]}

\nomenclature[P]{\( q \)}{Elementary electric charge}
\nomenclature[P]{\( \kb \)}{Boltzmann's constant}

\nomenclature[E]{\( V \)}{Electric voltage}
\nomenclature[E]{\( I \)}{Electric current}
\nomenclature[E]{\( P \)}{Electric power}
\nomenclature[E]{\( T \)}{Temperature}

\nomenclature[S]{pv}{Photovoltaic}
\nomenclature[S]{oc}{Open-circuit}
\nomenclature[S]{sc}{Short-circuit}
\nomenclature[S]{mp}{Maximum power}

\nomenclature[A]{dc}{Direct current}
\nomenclature[A]{PV}{Photovoltaic}
\nomenclature[A]{SDM}{Single-diode model}
\nomenclature[A]{SD}{Standard deviation}
\nomenclature[A]{NREL}{National Renewable Energy Laboratory}

\def\BibTeX{{\rm B\kern-.05em{\sc i\kern-.025em b}\kern-.08em
    T\kern-.1667em\lower.7ex\hbox{E}\kern-.125emX}}
\begin{document}

\title{
    A Practical Example of the Impact of Uncertainty on the One-Dimensional Single-Diode Model
}

\author{\IEEEauthorblockN{Carlos C\'ardenas-Bravo}
\IEEEauthorblockA{\textit{Univ. Grenoble Alpes} \\
\textit{Univ. Savoie Mont Blanc, CNRS, LAMA} \\
Chamb\'ery, France   \\
carlos.cardenas@cea.fr}
\and
\IEEEauthorblockN{Sylvain Lespinats}
\IEEEauthorblockA{
\textit{Univ. Grenoble Alpes} \\
\textit{CEA, Liten, INES} \\
Le Bourget du Lac, France \\
sylvain.lespinats@cea.fr}
\and
\IEEEauthorblockN{Denys Dutykh}
\IEEEauthorblockA{\textit{Mathematics Department} \\
\textit{Khalifa University}\\
Abu Dhabi, United Arab Emirates \\
denys.dutykh@ku.ac.ae}
}

\maketitle

\begin{abstract}
The state of health of solar photovoltaic (PV) systems is assessed by measuring the current-voltage (\iv) curves, which present a collection of three cardinal points: the short-circuit point, the open-circuit point, and the maximum power point. To understand the response of PV systems, the \iv{} curve is typically modeled using the well-known single-diode model (SDM), which involves five parameters. However, the SDM can be expressed as a function of one parameter when the information of the cardinal points is incorporated into the formulation. This paper presents a methodology to address the uncertainty of the cardinal points on the parameters of the single-diode model based on the mathematical theory. Utilizing the one-dimensional single-diode model as the basis, the study demonstrates that it is possible to include the uncertainty by solving a set of nonlinear equations. The results highlight the feasibility and effectiveness of this approach in accounting for uncertainties in the SDM parameters.
\end{abstract}

\begin{IEEEkeywords}
uncertainty, single-diode model (SDM), photovoltaics (PV), mathematical analysis
\end{IEEEkeywords}


\textbf{
© 2024 IEEE.  Personal use of this material is permitted.  Permission from IEEE must be obtained for all other uses, in any current or future media, including reprinting/republishing this material for advertising or promotional purposes, creating new collective works, for resale or redistribution to servers or lists, or reuse of any copyrighted component of this work in other works.
}

\section{Introduction}

\begin{figure}[b]
    \centering
    \includegraphics[width=2.5in]{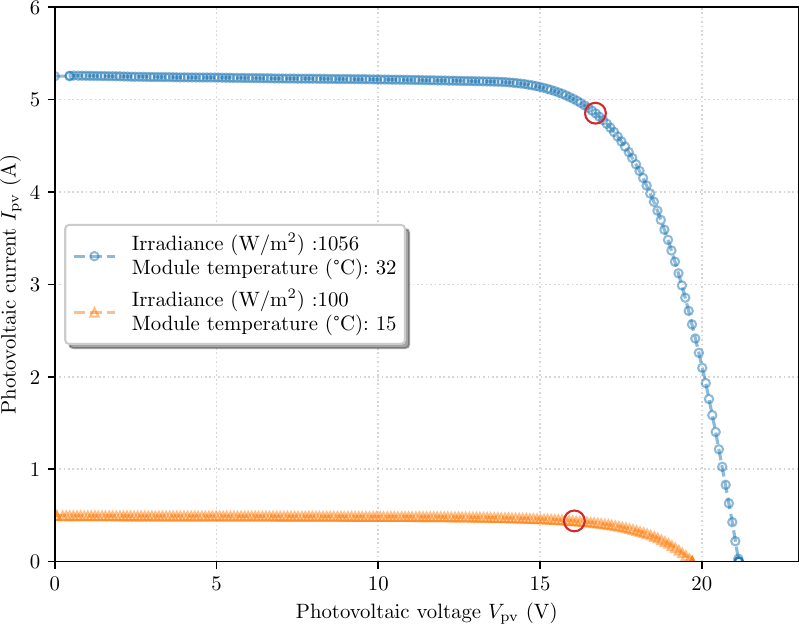}
    \caption{\iv{} curves samples representing two operational conditions.
    The maximum power point for each curve is highlighted in a red circle. 
    The data set used corresponds to the solar module Cocoa\_mSi460A8 \cite{marion_new_2014}.
    }
    \label{fig:iv-curve}
\end{figure}

In the field of solar photovoltaic (PV) energy, several diagnostic tools are available to assess the health status of PV systems at the module level. 
One of the most commonly used tools is the current-voltage (\iv) curve, which typically exhibits the shape shown in \refig{fig:iv-curve}. 
This curve is characterized by three key points named cardinal (or remarkable) points: (i) the short-circuit point \((0, \Isc)\), where \(\Isc\) is the short-circuit current; (ii) the open-circuit point \((\Voc, 0)\), where \(\Voc\) is the open-circuit voltage; and (iii) the maximum power point \((\Vmp, \Imp)\), where \(\Vmp\) is the voltage at maximum power and \(\Imp\) is the corresponding current. 
A thorough analysis of these cardinal points provides insights into the PV system's performance under varying irradiance and temperature conditions (see \refig{fig:iv-curve}). 
To perform this analysis, models are typically employed as tools to interpret the I-V curve.

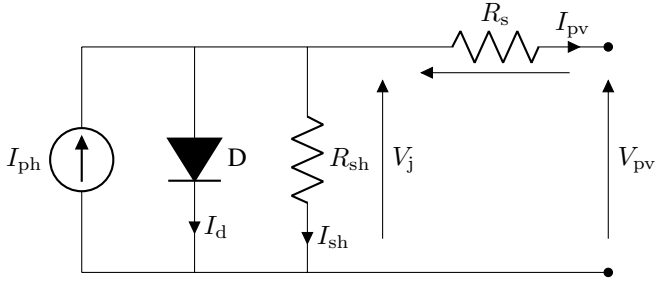
\begin{figure}[t]
    \centering
    \begin{circuitikz}[american,straight voltages]
    \draw 
        (0,0) to[isource, l=$\Iph$] (0,3)
        (1.5,3) to[full diode=D, i=$\Id$] (1.5,0)
        (3,3) to[R=$\Rshh$, i=$\Ish$] (3,0)
        (4,3) to[R=$\Rs$, -*, i=$\Ipv$, v<=$ $] (7,3)
        (4,0) to [open, v=$\Vj$] (4,3)
        (7,0) to [open, v=$\Vpv$] (7,3)
        (0,0) to[short,-*] (7,0)
        (0,3) to[short,-] (4,3)
    ;
    \end{circuitikz}
    \caption{ 
    Electrical diagram of the single-diode model. 
    Four electrical elements compose this model: a DC current source of magnitude $\Iph$, an ideal diode, D, a shunt resistance, $\Rshh$, and a series resistance, $\Rs$.
    } \label{fig:sdm_circuit}
\end{figure}

One of the most extended models is the single-diode model (SDM), which is described by four electrical elements \cite{townsend_method_1989}, as depicted in \refig{fig:sdm_circuit}.
There is one current source of magnitude $\Iph$, which represents the photo-generated current.
Next, there is an ideal diode, D, described by the Shockley equation \cite{shockley_theory_1949}
\begin{equation}
    \Id \defi \Io \parentheses{\exp{ \parentheses{\cfrac{q \, \Vj}{n \, \Ns \, \kb \, \Tpv \, } } } - 1 },
\end{equation}
where $\Io$ is called the dark saturation current of the diode; $n$ is the non-ideal factor of the diode; $\Ns$ is the number of solar cells connected in series; $\Tpv$ is the temperature of the PV system in kelvin; $\Vj$ represents the junction voltage across the diode in volts, following the electrical reference indicated in \refig{fig:sdm_circuit}; $q$ is the elementary electric charge equal to $1.602\,176\,634\e{-19}$ C; and $\kb$ is the Boltzmann's constant equal to $1.380\,649\e{-23}$ J/K.
For simplifying the mathematical notation, the equivalent factor of the diode, $A$, is introduced as
\begin{equation}
    A \defi \cfrac{n\, \Ns\,\kb\,\Tpv}{q}.
\end{equation}
Then, the third element corresponds to the shunt resistance, $\Rshh$, which represents the leakage current in the solar cell.
Often, the shunt resistance is expressed by the shunt conductance, $\Gsh$, defined as
\begin{equation}
    \Gsh \defi \cfrac{1}{\Rshh}.
\end{equation}
Finally, the series resistance, $\Rs$, corresponds to the resistive losses produced, for example, by the interconnections between two adjacent solar modules.
Mathematically, the SDM is expressed as
\begin{equation}
    \label{eq:sdm}
    \Iph = \Io \parentheses{\ej-1} + \Gsh \, \Vj + \Ipv,
\end{equation}
where $\Vj$ is equivalent to:
\begin{equation}
    \Vj \defi \Vpv+\Rs\,\Ipv.
\end{equation}
Equation \eqref{eq:sdm} shows that the SDM is dominated by five parameters: $\Iph$, $\Io$, $A$, $\Gsh$ and $\Rs$.

To identify the parameters of the SDM, optimization tools are employed.
For instance, some methodologies utilize the cardinal points to express the single-diode model as a function of one variable \cite{laudani_reduced-form_2013,toledo_geometric_2014,cardenas-bravo_parameters_2024}.
This representation is called one-parameter single-diode model, SDM-1, and an online platform, \textsc{PVmodel}, was developed to compute the SDM parameters following Toledo \etal{} approach \cite{galiano_pvmodel_2022}. 
One important factor to consider in these approaches is the uncertainty introduced by real measurements. 

In the domain of photovoltaics, three predominant sources of uncertainty are identified \cite{dirnberger_uncertainty_2010}: I-V curve measurement, irradiance levels, and array temperature. 
These uncertainties propagate through the measurement process, impacting the accuracy of both current and voltage measurements. 
As a result, rather than a singular optimal SDM, there exists a parameter space within which the SDM parameters lives. 
For modelers, quantifying the uncertainty in these parameters is essential, as it facilitates the precise computation of operational points.
The computation of uncertainty has been explored by authors such as R. Ben Messaoud \cite{ben_messaoud_extraction_2020}, who utilized the simulated annealing algorithm (SA) to estimate optimal SDM parameters and their associated uncertainties. However, that method relies on a metaheuristic algorithm, making the uncertainty computation demanding.

In this paper, we discuss the impact of uncertainty on the computation of the one-dimensional single-diode model. 
The performed analysis is static in the sense that one operational point of irradiance and temperature is analyzed. 
To test this, the open-source data set provided by NREL \cite{marion_new_2014} is utilized.
Specifically, we selected the multi-crystalline solar module Cocoa\_mSi460A8.
From this data set, the uncertainty percentages for the current and voltage are extracted from the theory coming from the one-dimensional single-diode model.
Therefore, it is not necessary the use of metaheuristic algorithm, making the computational process more effective.

The rest of the paper is organized as follows. Section \ref{section:one-dimensional_single-diode_model} presents the theoretical foundations of the one-dimensional single-diode model. Section \ref{section:data_set_description} describes the data used in this study. Section \ref{section:results_and_discussion} presents and discusses the results. Finally, Section \ref{section:conclusion} indicates the conclusions of this work.


\section{One-dimensional single-diode model} \label{section:one-dimensional_single-diode_model}

The one-dimensional single-diode model (SDM-1) is an implicit function derived when the cardinal points are incorporated into \eqref{eq:sdm}. Although it is not possible to obtain the optimal parameters of the SDM-1 from this representation, it is possible to compute the parameter domain such that all five parameters of the single-diode model are positive. For instance, if the SDM-1 is expressed as a function of the diode's equivalent factor, then the maximum value for the diode's equivalent factor is the minimum value between
\begin{equation}
    \parentheses{f_\mathrm{mp}} \big|_{R_\mathrm{s}=0} = 0,
\end{equation}
and the solution of the equation system
\begin{align}
    f_\mathrm{sh} = 0, \\
    f_\mathrm{mp} = 0,
\end{align}
where
\begin{multline}
    \label{eq:fsh}
    \fsh \defi \parentheses{\Isc-\Imp} \eoc \\- \Isc \emp + \Imp \esc,
\end{multline}
\begin{multline}
    \label{eq:fmp}
    f_\mathrm{mp} \defi -A  V_\mathrm{mp}  ( 2  I_\mathrm{mp} - I_\mathrm{sc})  \exp{\left( \cfrac{V_\mathrm{oc}}{A} \right)} \\ + \left( ( V_\mathrm{oc}  I_\mathrm{mp} + V_\mathrm{mp}  I_\mathrm{sc}  - V_\mathrm{oc}  I_\mathrm{sc} )  (V_\mathrm{mp} - R_\mathrm{s}  I_\mathrm{mp}) \right. \\ \left. +  A \, ( V_\mathrm{oc}  I_\mathrm{mp} - V_\mathrm{mp}  I_\mathrm{sc} ) \right)  \exp{\left( \cfrac{V_\mathrm{mp}+R_\mathrm{s} I_\mathrm{mp}}{A} \right)} \\ + A \, I_\mathrm{mp} \, ( 2  V_\mathrm{mp} - V_\mathrm{oc}) \, \exp{\left( \cfrac{R_\mathrm{s} I_\mathrm{sc}}{A} \right)}.
\end{multline}

\begin{figure}[b]
    \centering
    \includegraphics[width=2.5in]{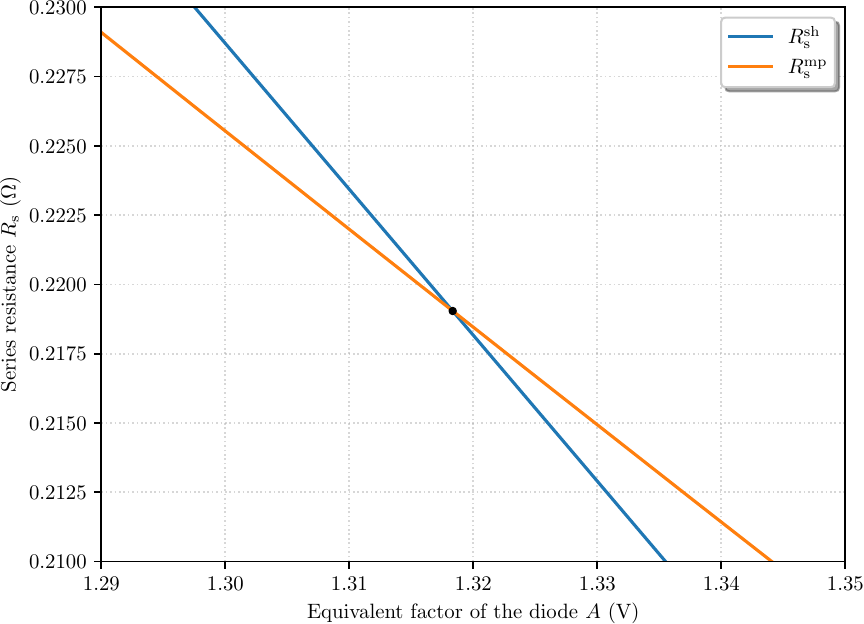}
    \caption{Maximum power series resistance, $\Rsmp$, and shunt series resistance, $\Rssh$ computed from the solar module Cocoa\_mSi460A8 (data reported in \cite{marion_new_2014}).
    }
    \label{fig:example-amax}
\end{figure}

All other possible solutions of the SDM-1 are derived from $\fmp$. 
More details about this formulation are indicated in \cite{cardenas-bravo_parameters_2024} with the codes hosted in the CodeOcean servers \cite{cardenas-bravo_supplementary_2024}. \refig{fig:example-amax} depicts the series resistance computed from \eqref{eq:fmp}, $\Rssh$, and the series resistance computed from \eqref{eq:fsh}, $\Rsmp$, both presented as functions of $A$. It can be seen that the functions intersect at the point $(1.3183\ \text{V},\, 0.2190\ \Omega)$. 
In turn, this point can be computed from different combinations of cardinal points. 
Therefore, it is possible to compute the effect of the uncertainty of the cardinal points on the parameter domain of the SDM-1.

\section{Data set description} \label{section:data_set_description}

The open-source data set provided by NREL \cite{marion_new_2014} is used for the current analysis.
This data set includes several PV technologies, such as multi-crystalline silicon (m-Si), cadmium telluride (CdTe), and copper indium gallium selenide (CIGS), among others.
For each technology described in the data set, there is a collection of I-V curves measured under real operating conditions. 
In addition, several pieces of information from the \iv{} curve are reported, such as the cardinal points and their uncertainty, time of measurement, irradiance, and cell temperature.
Specifically, the uncertainty is computed taking into account four components of uncertainty described as follows.
The first one is associated with the range and uncertainty of the multi-meter used to perform the measurement.
Next, the uncertainty of the cardinal point under analysis is included. 
Depending on the cardinal point, there is a different methodology on how to estimate it.
For instance, the uncertainty of the open-circuit voltage is estimated from the standard error of a linear fit in the \iv{} curve region such that $\Vpv>0.9\Voc$ and $\Ipv<0.2\Isc$.
The third component of uncertainty is the resolution of the I-V curve tracer.
Lastly, it has to be included the uncertainty of the standard cardinal point used in the calibration of the I-V curve tracer.

\begin{figure}[b]
    \centering
    \includegraphics[width=2.5in]{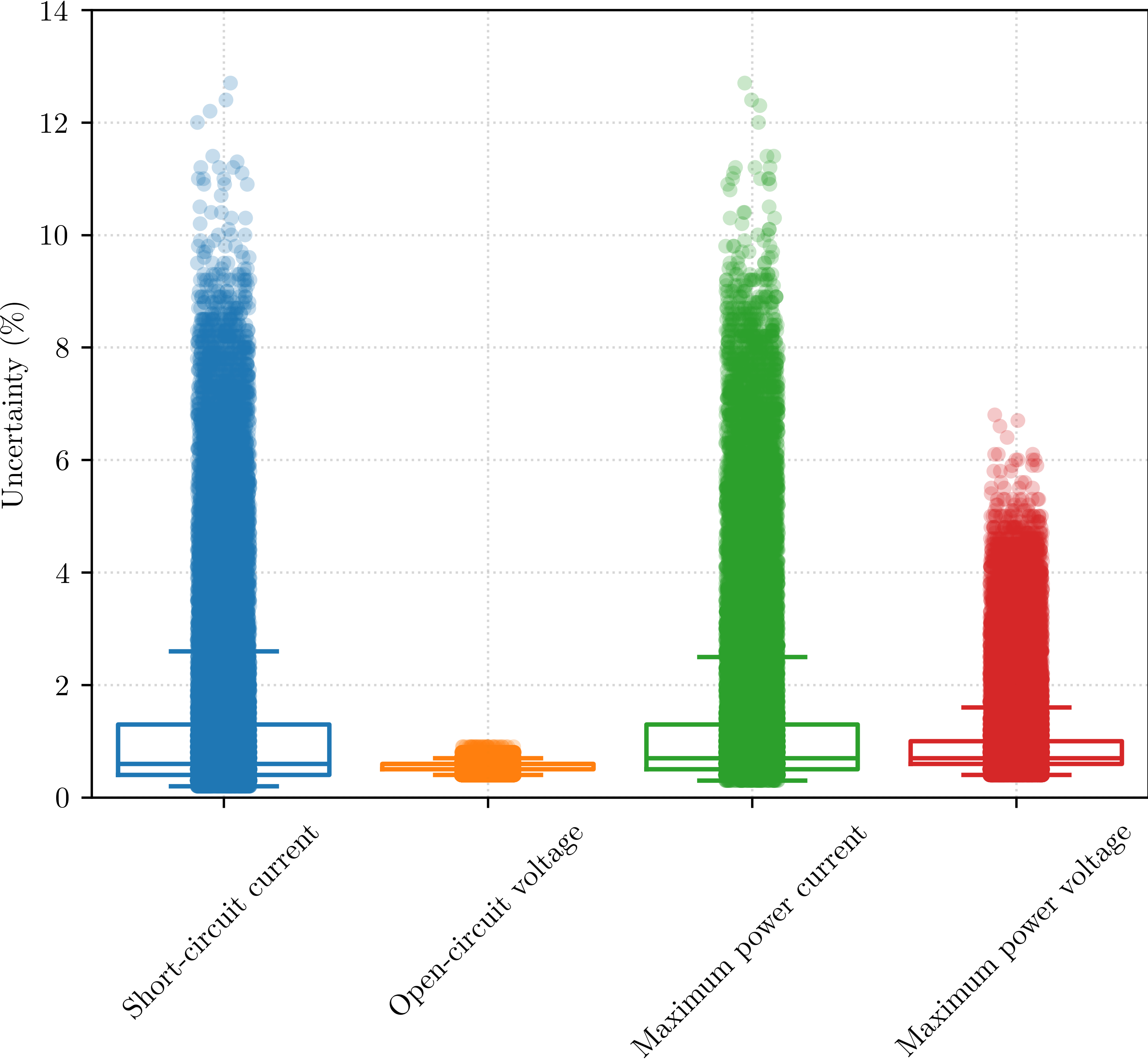}
    \caption{Box plot extracted from the uncertainties of the cardinal points. The data was extracted from \cite{marion_new_2014}.}
    \label{fig:dataset-boxplot-2024}
\end{figure}

From this database, the selected PV module corresponds to the m-Si Cocoa\_mSi460A8.
This choice is made since silicon technology is one of the most extended nowadays.
The data file of associated with this module contains 38\,929 I-V curves measured in Cocoa, Florida, from January 2021 to March 2012.
The uncertainties of the cardinal points on the data set are expressed as percentages, with a precision of one decimal place.
\refig{fig:dataset-boxplot-2024} depicts the boxplot computed from the uncertainties of the short-circuit current, the open-circuit voltage, the maximum power current, and the maximum power voltage.
From this figure, the biggest source of uncertainty is presented by the short-circuit current and the maximum power current, having a maximum value of 12.7\%.
Regarding the voltage, the maximum uncertainty is presented by the maximum power voltage, with a value of 6.8\%.

\retab{tab:dataset-summary-statistics} summarizes the minimum, mean, maximum, and standard deviation values computed from the cardinal points.
The highest uncertainty deviations corresponds to the current, with a value of 1.5\% for both the short-circuit current and the maximum power current.
It is noted that the uncertainty between the short-circuit current and the maximum power current is similar, due to the method described above.
On the other hand, the voltage presents the largest uncertainty deviation for the maximum power voltage, with a value of 0.8\%.

\begin{table}[t]
    \centering
    \caption{Summary statistics derived from the uncertainties of the cardinal points. The data was extracted from \cite{marion_new_2014}. } \label{tab:dataset-summary-statistics}
    \begin{tabular}{|c|c|c|c|c|}
    \hline
        Variable & Min & Mean & Max & SD  \\
    \hline
        Short-circuit current (\%) & 0.2 & 1.2 & 12.7 & 1.5 \\
    \hline
        Open-circuit voltage (\%)  & 0.4 & 0.5 & 0.9 & 0.1 \\
    \hline
        Maximum power current (\%) & 0.3 & 1.3 & 12.7 & 1.5   \\
    \hline
        Maximum power voltage (\%) & 0.4 & 1.0 & 6.8 & 0.8  \\
    \hline
    \end{tabular}
\end{table}

In this work, we illustrate the computation of uncertainty using the \iv{} curve measured on January 22, 2011, at 12:05:04 (see the blue circled curve in \refig{fig:iv-curve}.) 
The cardinal points and their corresponding uncertainties (measured in A or V, as appropriate) are summarized in \retab{tab:summary-cardinal-points}. 
Detailed information about the \iv{} curve can be found in \cite{marion_users_2014}.

\begin{table}[t]
    \centering
    \caption{Cardinal points and its uncertainties for the \iv{} curve measured on January 22, 2011, at 12:05:04. 
    The solar module Cocoa\_mSi460A8 is used as a basis and the data is reported in \cite{marion_new_2014}. } \label{tab:summary-cardinal-points}
    \begin{tabular}{|c|c|c|}
    \hline
        Variable & Value    \\
    \hline
        Short-circuit current (A) & $5.26\pm0.02$ \\
    \hline
        Open-circuit voltage (V)  & $21.15\pm0.08$ \\
    \hline
        Maximum power current (A) & $4.85\pm0.02$   \\
    \hline
        Maximum power voltage (V) & $16.71\pm0.07$  \\
    \hline
    \end{tabular}
\end{table}

\section{Results and Discussion} \label{section:results_and_discussion}

\begin{figure}[b]
    \centering
    \includegraphics[width=2.5in]{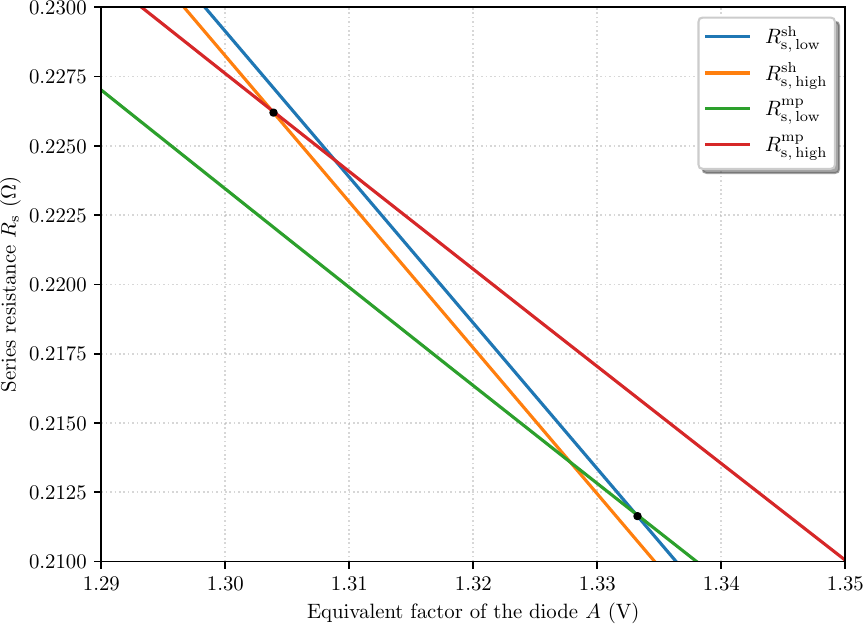}
    \caption{Maximum power series resistance, $\Rsmp$, and shunt series resistance, $\Rssh$ computed for the low and high \iv{} curves.
    The solar module Cocoa\_mSi460A8 is used as a basis and the data is reported in \cite{marion_new_2014}.
    }
    \label{fig:uncertainty-region}
\end{figure}

From the cardinal points presented in \retab{tab:summary-cardinal-points}, two pairs of functions, $\Rssh$ and $\Rsmp$, can be derived. The first pair of functions, $\Rsshlow$ and $\Rsmplow$, is computed from the lowest collection of cardinal points. Conversely, the second pair of functions, $\Rsshhigh$, and $\Rsmphigh$, is computed from the highest cardinal points. Each pair of points yields a specific combination of $(\Amax,\,\Rsmin)$. \retab{tab:summary-low-high} summarizes the low and high cardinal point values used for computing the two pairs of functions, $\Rssh$ and $\Rsmp$. \refig{fig:uncertainty-region} illustrates these functions.
Here, the intersection point for the low case is $(1.3332,\,0.2116)$ and the for the high case is $(1.3039,\,0.2262)$.
This fact suggest that there is a region where the parameters $\Amax$ and $\Rsmin$ belong.

\begin{table}[t]
    \centering
    \caption{Summary of cardinal points for low and high \iv{} curves. 
    The solar module Cocoa\_mSi460A8 is used as a basis and the data is reported in \cite{marion_new_2014}. } \label{tab:summary-low-high}
    \begin{tabular}{|c|c|c|}
    \hline
        Variable & Low & High \\
    \hline
        Short-circuit current (A) & 5.25 & 5.28 \\
    \hline
        Open-circuit voltage (V)  & 21.06 & 16.78 \\
    \hline
        Maximum power current (A) & 4.83 & 4.87 \\
    \hline
        Maximum power voltage (V) & 16.65 & 16.78  \\
    \hline
    \end{tabular}
\end{table}

\begin{figure}[t]
    \centering
    \includegraphics[width=2.5in]{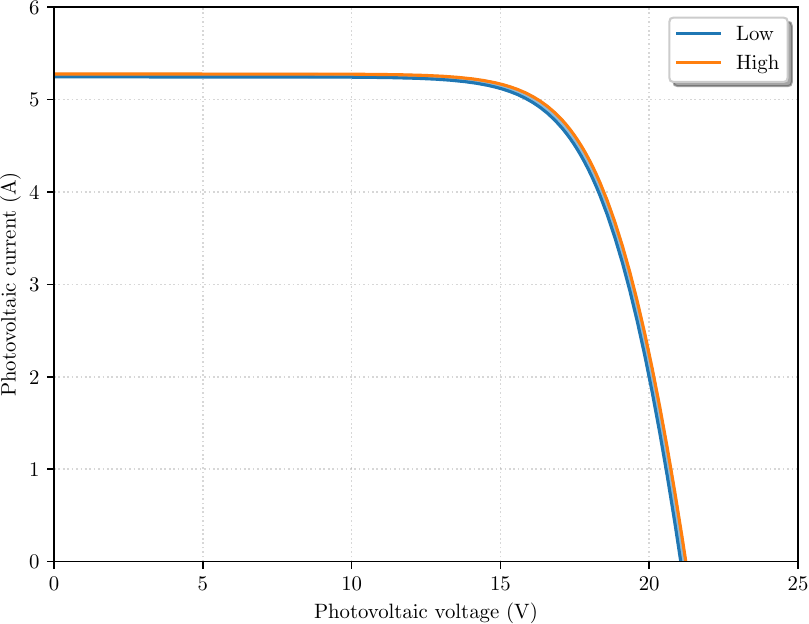}
    \caption{
        \iv{} curves computed for the low and high combinations of cardinal points.
    }
    \label{fig:iv-curves}
\end{figure}

The impact on the \iv{} curves can be visualized using the aforementioned parameter combinations $(\Amax,\,\Rsmin)$. 
This is shown in \refig{fig:iv-curves}, highlighted in blue. 
Consequently, the uncertainty in the cardinal points necessarily implies that the maximum value of the equivalent factor of the diode belongs to a closed interval.
For this example, the interval of $\Amax$ corresponds to $[1.3039,\,1.3332]$. 
Regarding the other four parameters $\Iph$, $\Io$, $\Rs$, and $\Gsh$, since they are functions of $A$, it is not ensured that the uncertainty is a fixed value.
Instead, it could correspond to a function parameterized with $A$.

\section{Conclusions} \label{section:conclusion}

This paper presents the computation of the uncertainty in the one-dimensional single-diode model (SDM-1).
The SDM-1 is expressed as an implicit function of the equivalent factor of the diode, \( A \). 
An illustrative example using the Cocoa\_mSi460A8 module, extracted from the National Renewable Energy Laboratory (NREL) dataset, is provided. 
The dataset reports the uncertainty of the cardinal points. 
The results for this specific example indicate that the equivalent factor of the diode belongs to a defined region. 
The findings suggest that the uncertainty of the other parameters, $\Iph$, $\Io$, $\Rs$, and \(\Gsh\), depends on the \( A \) parameter. 
Practical applications for the uncertainty region include enhancing the accuracy of operational points and conducting fault detection analysis based on the parameter uncertainties.

To determine these regions with greater accuracy, future work should focus on calculating multiple points \((\Amax,\,\Rsmin)\) derived from the intersection of the functions \(\Rssh\) and \(\Rsmp\). This approach will offer a more comprehensive understanding of the uncertainties in the SDM parameters. Additionally, future research should explore alternative models for representing the \iv{} curve, such as the double-diode model or the single-diode model expressed in terms of the five original parameters. Incorporating multiple irradiance levels and temperatures into these models will further enhance the robustness of the analysis.

\section*{Acknowledgment}

The work of D. Dutykh has been supported by the French National Research Agency, through the Investments for Future Program (ref. ANR-18-EURE-0016 - Solar Academy).  In addition, C. C\'ardenas-Bravo would like to acknowledge the INES.2S French Institute for the Energy Transition, the PhD's contract USMB/2021-454, and the CSMB --- Conseil Savoie Mont Blanc for its support. This publication is based upon work supported by the Khalifa University of Science and Technology under Award No. FSU-2023-014.

\bibliographystyle{IEEEtran}
\bibliography{IEEEabrv,references}

\end{document}